# 'Congruent Partitions' of Polygons – a Short Introduction

-R. Nandakumar
(nandacumar@gmail.com)

**Abstract:** We introduce the problem of partitioning 2D regions (usually convex polygonal regions) into mutually congruent pieces.

## 1. Introducing the Problem

This article is an elaboration on [1] where this problem was originally stated. It is essentially a collection of inter-related questions and a few partial answers.

We recall a definition from basic Euclidean geometry:

Two planar regions are **congruent** if one can be made to perfectly coincide with the other by translation, rotation or reflection (in 2D, flipping over).

**The Basic Problem:** To partition a given polygonal region *P* into N mutually congruent pieces (or 'tiles') so that the fraction of the area of *P* not covered by the union of the pieces is as small as possible. N is a finite positive integer. Each tile should have finitely many sides.

**Some Definitions:** A partition which leaves out the least area from P is an *optimal* congruent partition for that N. If a congruent partition leaves no area of P unused, then it is called a *perfect congruent partition*.

If for a given input polygon and N, the optimal congruent partition is not unique, we are interested in finding the one in which the pieces have least complexity (least number of edges) – such a partition may be called the *simplest* optimal congruent partition.

*Remark:* In the following discussions, the region to be partitioned P and the tiles are convex, unless mentioned otherwise.

## 2. Preliminary Results and Questions

***1. It is known that there exist even very simple polygons which do not allow perfect congruent partitioning (no left over) into pieces of finite complexity for any N.***

Reference [2] proves that a convex quadrilateral with angles linearly independent over the rationals (for example, the 4 angles could be, in degrees, {$\alpha_1=180/\sqrt{5}$, $\alpha_2=180/\sqrt{7}$, $\alpha_3=180/\sqrt{11}$, $\alpha_4=360-\alpha_1-\alpha_2-\alpha_3$ }) does not allow a perfect congruent partition for any N.

*Remark:* [2] actually tries to prove a somewhat stronger result: even pieces with same *sets of values* for their angles (the sequence of the angles in different pieces could be different, as could be their sizes) fail to fully fill out the specified quadrilateral.

***2. It is almost certain that if P convex, the simplest tile that gives the optimal congruent partition for some N need not be convex as well.***

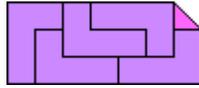

Consider the above example due to Friedman [3]. From a 3x7 rectangle, we chop off an isosceles triangle of area 1/2 from one corner, leaving a convex pentagon of area 20.5. For N=5, this pentagon allows a layout of 5 L-tetrominos (a non-convex shape of area 4) which leaves out only 0.5 units of area. There appears to be no 'better' congruent partitioning for N=5 with convex tiles.

***3. A Claim:*** Given a convex P and any N; if P allows a perfect congruent partition of itself into N *non-convex* pieces each with finitely many sides, then P also allows a perfect congruent partition into N convex pieces with finitely many sides.

This claim holds for N=2. For N=3 and beyond, things are not clear.

A convex region, which allows perfect congruent partition *only* with non-convex tiles with finitely many sides for some specific N would counter this claim. On the other hand, if the above claim is true in general, that would mark a qualitative difference between perfect and non-perfect congruent partitions. It also puts strong constraints on the complexity of the simplest piece which achieves a perfect congruent partition of any convex polygon (see **appendix** below) and hence could help design a practical algorithm to check for the existence of perfect congruent partitions for a given convex polygon. *Question:* Among all perfect congruent partitions of any given convex P, is the partition that minimizes the sum of perimeters of pieces, a partition into convex tiles?

***4. Upper-bounding the area of the input P that gets left-over:*** The most basic question in this direction could be the following:
*What is the shape of a convex region (with boundary not necessarily formed by straight edges) P such that if it is optimally congruent partitioned into 2 convex pieces, the fraction of the area of P that is left over is a maximum?*

In other words, If P and the tiles are convex, can there be bounds on how sub-optimal a convex congruent partition into 2 pieces (in general, N pieces) can be?

***5. Non-convex tiles:*** Question: if the tiles are allowed to be non-convex and *arbitrarily complex*, can we always achieve, for any N, a 'near-perfect' congruent partition – a congruent partition that is arbitrarily close to being a perfect congruent partition? Such a partition makes the fraction of P's area that goes waste *tend to 0*. For instance, the N tiles could be tightly packed spirals radiating from a suitable core region point in the interior; the 'pitch' of each spiral going to zero (causing the spiral to wind around the 'focal region' infinitely many times).

***6.* The 3D version of the problem:** If we have 3D regions rather than polygons to be partitioned, is it true that *any* 3D region can indeed be partitioned into N mutually congruent connected regions for *any* N with no leftover *or the leftover tending to 0* – with the pieces arbitrarily complex and also be densely entangled with one another (a similar scenario is discussed in [4], where an apple is shown to be eaten by 2 species of thin worms with the species avoiding one another – resulting in 2 separate connected regions of infinite complexity). In 2D (item 5 above) where there is no scope for entanglement of pieces, such a partition would be more difficult. Indeed, *entanglement* of regions could be a feature that could set perfect congruent partitions in higher dimensions apart from 2D.

We could also consider the 3D version of the claim in item 3 above. Consider a 3D region R that does not allow perfect congruent partition into N *convex* 3D regions. As suggested above, one may well be able to use entanglement to achieve a congruent partition of it into 'infinitely fibrous' pieces that is perfect. But is there such an R which cannot have perfect congruent partitions into N convex 3D regions but can be perfect congruent partitioned into N non-convex pieces of *finite complexity*?

***7. Deciding whether a given polygon P (not necessarily convex) allows a perfect congruent partition into N pieces of finite complexity*** is another variant of our problem. The case N=2 has been solved in [4]. For larger N, things appear uncertain, even for convex P.

*A possibility:* If a convex polygon allows perfect congruent partition into N convex pieces, then a layout could always be created with only a few types of basic topology for any N – a straight chain, star, etc.. If that is indeed the case, the decision of whether a given convex polygon allows perfect congruent partition with N *convex* pieces, can be achieved by considering a few simple and enumerable cases.

*8.* Rather than **break** the input polygon P into congruent pieces, we could try to find a piece such that N congruent copies of it **cover** the whole of P with the least area of the *tiles* going waste (parts of the tiles could extend beyond of the boundary of P or may be allowed to **overlap** among themselves).

*9.* Instead of a set of N mutually congruent pieces, we could ask about partitioning a given region into N mutually identical *sets of pieces* (with each set required to have some finite number of pieces).

A simple fact: any triangle can be cut into 4 mutually congruent triangles. And any polygon can be divided into triangles; any m-gon will give m-2 triangles.
This implies: if N is 4, we can triangulate *any* (m-vertex) input polygon P and further divide each resulting triangle into 4 so as to achieve a partition into 4 identical tile-sets each with m-2 triangles. This is easily generalizable to N = any perfect square. No bit of P goes waste.

*10.* If we restrict the definition of 'congruence' of the pieces to only translation and rotation (this will mean a polygon and its mirror image are not necessarily congruent), how will the problem change?

# 3. Acknowledgements

The author is grateful to Swami Sarvottamananda (Shreesh Mj) and Br. Swathy Prabhu for discussions.

# 4. References

1. http://maven.smith.edu/~orourke/TOPP/P73.html (May 2009)
2. http://domino.research.ibm.com/Comm/wwwr_ponder.nsf/challenges/December2003.html (December 2003)
3. Erich Friedman – personal communication (May 2009).
4. George Gamow - "One, Two Three... Infinity"
5. Dania El-Khechen, Thomas Fevens, John Iacono, and Günter Rote.
   Partitioning a polygon into two mirror congruent pieces.
   In *Proc. 20th Canad. Conf. Comput. Geom.*, pages 131-134, August 2008

# Appendix

Here we derive a simple result on the complexity (number of edges) of tiles; the arguments follow [2].

**Claim:** If a given convex polygonal polygonal region P allows a perfect congruent partition for any given n with n *convex* tiles, the complexity of the tile is limited by that of P itself, with no dependence on n.

**Proof**:

Some notation:
    $p$ – the number of vertices on the full convex polygon P.
    $n$- the number of mutually congruent convex tiles which perfectly partition the full polygon.
    $k$- the number of vertices on each tile.
    $r$- the number of vertices in the layout of tiles which lie on the boundary of the full polygon P, which are not the vertices of P itself (call these the boundary vertices in the layout).
    $m$- the number of vertices in the layout which lie in the interior of the full polygon (call these the internal vertices in the layout)

We obtain the following relations:

$$nk >= 3m + 2r + p \quad (1)$$

$$n(k-2) = 2m + r + p - 2 \quad (2)$$

Note: Relation (1) follows from the fact that when tiles are convex, every internal vertex in the layout has at least 3 tiles meeting there; further, at every boundary vertex, at least 2 tiles meet. *nk* is the total number of vertices for all the tiles taken together. Equation (2) follows from equating two sums of the angles, measured in units of π; at each internal vertex in the layout the sum of angles meeting is 2 π and at every boundary vertex, π .

(1) and (2) together yield the inequality:

$$(6 - k)n >= r - p + 6. \quad (3)$$

We now consider the possibility of the complexity of the tile - its number of edges *k* -, being larger than the complexity of the full polygon (p).
So, we set $k = p + \alpha$ where $\alpha$ is a +ve integer.
With this substitution, (3) yields, $6 - p - \alpha >= r/n + 6/n - p/n$.

$$=> p(n-1)/n <= 6 - 6/n - \alpha - r/n. \quad (4)$$

Now, the expression on the right side of (4) cannot be greater than 5 (because, by definition, $\alpha >= 1$). Since $(n-1)/n$ cannot be less than than ½ for positive integer *n*, we get $p <= 10$ for any value of n.

Thus, if the convex tile is to have exactly 1 edge more than P, the complexity of P and hence, the tile are limited by 10. And if the tile is to have a greater excess of edges $\alpha$, the right side of (4) has a still lower limit than 5 so *p* and *k* are limited to 8 and 10 respectively.

*Conclusion:* if a tile in a convex congruent partition is to be more complex than a convex P, neither of them can be very complex. Of course, if α can be negative (if P is allowed to be more complex than the tile), *k* can be arbitrarily large - but less than *p,* the complexity of P.

*Note:* for very small values of *p* and *n*, we sometimes have tiles more complex than P itself: for example, if P is an equilateral triangle and *n* = 3, we can have a perfect congruent partition with 3 quadrilaterals *(k=4* and *p=3).*